%% file: DLS-MIP.tex
\documentclass[12pt,authoryear,times,5p,twocolumn]{elsarticle}
\usepackage{graphicx, booktabs, bigstrut, multirow,rotating,tcolorbox} 
\usepackage{amsmath,amsfonts,amssymb,amsthm} 
\usepackage[T1]{fontenc} 
\usepackage[english]{babel} 
\usepackage[export]{adjustbox} 
\usepackage{hyperref}

\newcommand{\set}[1]{\left\{#1\right\}}
\newcommand{\setn}{\{1,2,\dots,n\}}
\newcommand{\sett}[1]{\{1,2,\dots,#1\}}
\newcommand{\meq} {\stackrel{!}{=}}

\def\real{\mathbb{R}}
\def\binary{\mathbb{B}}
\def\setN{\mathbb{N}}

\newcommand{\smat}[1]{\ensuremath{\left[ \begin{smallmatrix} #1 \end{smallmatrix} \right]}}
\newcommand{\tp}[1]{\ensuremath{#1^\intercal}} 
\newcommand{\lequ}[2]{\begin{equation}#2\label{eq:#1}\end{equation}} 
\newcommand{\obj}[2][max]{& \text{#1} & & #2 \\}  
\newcommand{\objd}[3][max]{&\underset{#2}{\text{#1}} & & #3 \\}  
\newcommand{\sto}[1]{& & & #1 \\}  
\newcommand{\LPs} [2]{\begin{aligned} #1 & \text{subject to} & & \\ #2 \end{aligned}}
\newcommand{\LP}  [2]{\begin{equation}\LPs{#1}{#2} \end{equation}}
\newcommand{\mcomment}[1]{\text{\textcolor{cyan}{~~$\triangleright$\textit{~#1}}}} 
\newcommand{\ie}{\Rightarrow}
\newcommand{\equi}{\Leftrightarrow}
\newcommand{\ones}[1][n]{\tp{(1,1,\dots,#1)}}
\newcommand{\qbox}[1]{\begin{center}\fcolorbox{ocre}{ocre!10}{#1}\end{center}}
\newcommand{\qtable}[5]{
	\begin{table}[htbp]	\centering 
		\caption{#1}	\label{tab:#2}
		\begin{tabular}{#3}\toprule #4 \\	\midrule #5	\bottomrule	\end{tabular}
	\end{table}
}
\definecolor{ocre}{RGB}{243,102,25}
\newcommand{\tcb}[1]{\textcolor{blue}{#1}}

\theoremstyle{definition}
\newtheorem{example}{Example}[section]

\journal{arXiv}
\title{Dynamic lot size MIPs for multiple products and ELSPs with shortages, capacity and changeover limits}
\author{Wolfgang Garn} \date{April 2020}
\address{The Surrey Business School, University of Surrey, Guildford, Surrey, GU2 7XH, United Kingdom}
\ead{w.garn@surrey.ac.uk} \tnotetext[phone_and_fax_number]{\\ \\Tel.: +44(0)1483 68 2005; fax: +44(0)1483 68 9511.}

\begin{document}
	\begin{abstract}
Scheduling multiple products with limited resources and varying demands remain a critical challenge for many industries. 
This work presents mixed integer programs (MIPs) that solve the Economic Lot Sizing Problem (ELSP) and other Dynamic Lot-Sizing (DLS) models with multiple items. 
DLS systems are classified, extended and formulated as MIPs.
Especially, logical constraints are a key ingredient in succeeding in this endeavour. 
They were used to formulate the setup/changeover of items in the production line. 
Minimising the holding, shortage and setup costs is the primary objective for ELSPs. 
This is achieved by finding an optimal production schedule taking into account the limited manufacturing capacity. 
Case studies for a production plants are used to demonstrate the functionality of the MIPs. 
Optimal DLS and ELSP solutions are given for a set of test-instances. 
Insights into the runtime and solution quality are given. 
\end{abstract}

\maketitle

\section{Introduction}
Agile manufacturing is essential to industries where demand varies.
Dynamic Lot Sizing (DLS) models deliver production schedules that accommodate dynamic time dependent demand.
The aim of this paper is to present a general DLS MIP formulation.

This section develops a taxonomy of Dynamic Lot-Sizing (DLS) models 
and establishes the context to Economic Lot Sizing Problems (ELSP).
Used solution methods will be briefly mentioned with a focus on MIP approaches.
Essential insights and problem formulations are extracted using a consistent notation throughout this paper.
Occasionally, new formulation DLS MIP formulations will be given.

\subsection{DLS models}
In 1958 \citet{wagner2004dynamic} produced a ``landmark'' paper introducing a dynamic version of the 
Economic Lot Size (or Economic Order Quantity) model. 
Varying deterministic demand $d_j$ has to be fulfilled to create the production schedule $q=(q_1,\ldots,q_n)$ for one product over $n$ periods. 
The inventory entering period $t$ is $I=I_{0}+\sum_{j=1}^{t-1} q_{j}-\sum_{j=1}^{t-1} d_{j} $,
where $I_0$ is the initial stock, $q_t$ the manufactured goods and $d_t$ the demand during period $t$.
They formulated a minimal cost policy with holding costs $c^h_t$ and setup costs $c^o_t$ for periods $t \in \set{1,\dots,n}$: 
\lequ{WagnerBell}{f_{t}(I)=\min _{q_{t} \geq 0 \atop I+q_t \geq d_t} \left[c^h_{t-1} I+\omega_t c^o_t+f_{t+1}\left(I+q_t-d_t\right)\right],}
where $\omega_t=1$ represents setup occurring when $q_t>0$ and zero otherwise.
As can be seen the model does not allow shortages and is recursive.
They postulated three basic theorems about characteristics of an optimal program.
These allowed them to rewrite the policy as an more efficient recursive minimal cost program.
I would like to give an alternative MIP formulation and call it
\textit{ Dynamic Lot Size }(DLS) model \textit{with one product} (\textit{DLS-1p}):
\LP{ \label{LP:DLS}\objd[min]{q,\omega}{\tp{I} c^h  +\tp{\omega} c^o}}
   { \sto{I_t=I_{t-1}+q_t-d_t,~ t \in T = \sett{n},}
   	 \sto{\omega\leq Mq, ~\omega \geq mq, }
   	 \sto{q\in \real^n_+, \omega\in\binary^n,I_0\in\real_+,}
   }
where $e=\ones$, $M=\tp{e}d$ and $m=\frac{1}{M}$ ($M$ and $m$ can be tightened).
We used $\omega\leq Mq, ~\omega \geq mq$ to represent the logical constraints:
\qbox{\tcb{If} $q_t>0$ \tcb{then} $\omega_t = 1$ \tcb{else} $\omega_t=0$.}
Alternatively, the constraints 
\lequ{bigM}{q_t\leq M\omega_t,~t\in T} 
could have been used.
These are valid because of the following arguments.
If $q_t=0$ then $\omega_t=0$ or 1, 
but due to the minimisation objective $\omega_t=0$.
If $q_t>0$ then $\omega_t \meq 1$, because of (\ref{eq:bigM}).
If we assume that $\omega_t=0$ then (\ref{eq:bigM}) ensures that $q_t=0$.
On the other hand, if $\omega_t=1$ then $q_t<M$.
This proves the validity of (\ref{eq:bigM}) resembling the logical constraint.

\begin{example}{Wagner and Within.}
Let $d=\tp{\smat{69&29&36&61&61&26&34&67&45&67&79&56}}$, 
$c^h=\tp{\smat{1&1&\ldots&1}}$ 
and $c^o=\tp{\smat{85&102&102&101&98&114&105&86&119&110&98&114}}$.
This example was used in 1958 to illustrate the functionality of the recursive minimal cost program.
It gave the optimal cost 864. The here introduced MIP formulation gives the same result.

The special steady state data $d = 52.5e$, $c^h=e$ and $c^o=102.8e$, 
i.e. $\bar{d}=\frac{1}{n}\tp{e}d = 52.5$, $\bar{c}^h = 1$ $\bar{c}^o=102.8$
returns the optimal production quantity 105 every other week.
This is comparable to the economic order quantity (EOQ)
$Q=\sqrt{2\bar{d}\bar{c}^o/\bar{c}^h}\approx 103.9$, 
which is valid in the continuous case (i.e. the period $t\rightarrow0$).
\end{example}

The DLS-1p can be extended to a \textit{DLS with multiple products (DLS-mp)}.
\citet{boctor2004models} review models and algorithms for the \textit{dynamic-demand joint replenishment problem (DJRP)}.
This is a DLS-mp with an additional ``joint'' ordering cost (or common setup cost) regardless of the product(s) that are ordered (produced). 
The common ordering cost are known as general/common setup or changeover cost in a manufacturing context. 
This model will be abbreviated with \textit{DLS-mp-cs} instead of DJRP to be consistent with previous notations. 
The MIP after converting it to this paper's notation and using logical IP constraints is:
\LP{ \objd[min]{q,\Omega, \omega}{
		\sum_{t \in T}\left[C^o_{t} \Omega_{t}+\sum_{p \in P} \left[c^o_{tp} \omega_{tp}+c^h_{tp} I_{tp}\right]\right]}}
{   \sto{I_t=I_{t-1,p}+q_{tp}-d_{tp},~ t \in T,p\in P,}
	\sto{q_{tp}\leq M\omega_{tp}, ~\sum_{p\in P} \omega_{tp} \leq m \Omega_t,}
	\sto{Q=(q_{tp})\in \real^{n\times m}_+, (I_{tp}) \in \real^{n\times m}, I_0\in \real^p_+,}
    \sto{(\omega_{tp}) \in \binary^{n\times m}, \Omega\in \binary^p_+.}}
Here, a common setup (ordering) cost $C^o_t$ is used 
in addition to the individual setup (ordering) cost $c^o_{tp}$
for product $p\in P=\sett{m}$
with corresponding decision variables $\Omega_{t}$ and $\omega_{tp}$.
This formulation has two interesting features: (1) the common ordering cost,
and (2) the logical constraints.
When the common ordering cost is zero the program (\ref{LP:DLS}) defines the DLS-mp.
However, the number of decision variables must be reduced by dropping $\Omega$
and the following constraints have to be used: 
\lequ{limiting}{\omega_{tp} \leq M q_{tp}, ~\omega_{tp} \geq m q_{tp}.}
By the way, the objective function can be rewritten using matrices/vectors 
(for readability transpositions are implicit):
\lequ{obj-DLS-mp-cs}{C^o\Omega+\sum_{p \in P}(c^o_p\omega_p+c^h_p I_p).}

Next, I will mention a few works that illustrate other DLS-mp flavours.
\cite{gilbert2000coordination} looked at creating production schedules for
multiple products with a single production line. 
A production capacity requirement $\sum_{p\in P}q_{tp} \meq \hat{q}_t$ was introduced. 
Moreover, the model introduced constant priced products $(r_p)$, 
which allowed the formulation as a profit maximisation program.
Their revenue component defines a revenue/demand parameters $\beta_{tp} := d_{tp}/d_p$,
where $d_p$ was used for a demand intensity (related to price). 
This makes the price demand dependent.
Price changes lead to demand changes due to price elasticity.
However, it was assumed that this will not happen.
Moreover, the program gets the costs by solving the DLS with common setup costs.
Overall, some interesting first steps into considering price adaptations
and more importantly (in regard to DLS) a capacity requirement.
The capacity requirement will be picked up in Section \ref{sec:MIP} as production limit.

\cite{bruggemann2000discrete} introduced a DLS with batch production requirements.
This addresses an essential short-coming of the previous models.
Their model formulation assumes that a batch is produced
over several consecutive periods. 
This requires the introduction of decision variables and constraints
indicating the start and end of the production period.
The increased complexity motivated them to introduce
a ``two-phase'' (1. feasible solution search, 2. cost optimisation) heuristic based on simulated annealing.
Test-instances (products,periods)$ = \set{(3,30),(6,60),(8,80),(10,100)}$ 
were used to gain performance insights. 
The (3,30) instance was compared to the optimal solution having a 12.53\% difference.
I will propose a different approach to consider batch-sizes. 
Rather than using a small period (e.g. hours) 
that spans overs several periods, 
a large period (e.g. week) is used that allows the production of several batches.
This approach will be defined later in more detail
and overcomes the issue of the large number of decision variables.
Another idea to enhance their work are ``batch-interruption'' penalties in the objective function.
Overall, such DLS models are classified as \textit{DLS with batch production} (\textit{DLS-bp}) systems.

\cite{absi2009multi} proposed a multi-item capacitated lot-sizing problem with safety stocks and demand shortage costs.
They used a Lagrangian relaxation algorithm and a dynamic programming algorithm to solve this problem.
Generally, \textit{DLS models with capacity limit} will be classified as \textit{DLS-cl} systems.
Safety stocks are essential and their immediate consideration
when creating the production schedule seems promising.
This addresses shortages indirectly.
Interestingly, most DLS models found in this review do not allow shortages.
An alternative is to introduce safety buffers after having created the production schedule.
Here, the uncertainty (e.g. $N(\bar{d}_p,\sigma_{dp})$) in demand can be taken into account for each product
and the uncertainty in lead time (e.g. $N(\bar{L}_p,\sigma_{Lp})$).
Hence, the safety buffer is: 
\lequ{safetyBuffer}{b_p = Z\sqrt{\bar{d}_p^2 \sigma_{Lp}^2 +  \bar{L}_p^2 \sigma_{dp}^2}.}
Please see \cite[p223]{Garn2018Dec}, \cite[p121ff]{ghiani2004introduction} or \cite{silver1998inventory} for more details.

\citet{lu2011dynamic} provided a dynamic lot sizing for multiple products introducing a joint replenishment model,
where production quantities are of same size for all products per periods $q_{ti} = q_{tj}, i,j\in P$, 
i.e. $q_t$ can be used. 
Lost sales are considered by rejecting production requirements, i.e. no back-ordering.
This implies that no shortage decision variables are required for this formulation.
However, the rejected quantity $y_{tp}$ has to be recorded.
In addition to a common setup cost a common production cost $\pi_t$ is used.
The resulting MIP \footnote{The transpose operator has been used implicitly to improve the readability, i.e. $C^o \Omega := \tp{(C^o)} \Omega = \sum_{t\in T} C^o_t \Omega_t $, similarly variables non-negativity was assumed} is:
\LP{ \obj{C^o\Omega + \pi q_t + \sum_{p \in P} (c^h_p I_p + c^s_p y_p) }}
   { \sto{ I_{tp} = I_{t-1,p}+y_{tp} + q_t - d_{tp},} 
   	 \sto{q\leq M \Omega,~Y\leq D.} }
NP-hardness of the above problem was shown using the subset-sum approach.
The NP-hardness motivated the development of a heuristic.
The \textit{Dominant Product Approach} was used, 
which is actuated from the real-world observation 
that a few products generate the majority of economic benefits.
This is commonly known as the Pareto effect.
The proposed algorithm chooses the dominant product to generate the initial schedule.
Roughly speaking, the production of the remaining products is mainly guided by rejecting quantities.
The second heuristic aggregates demand of all products during a time period; and runs a LP. 

\citet{lee2005heuristic} gave a heuristic algorithm for a multi-product dynamic lot-sizing and shipping problem.
Basically, the meaning of the common ordering cost $C^o$ was substituted by freight cost of a container.
The difference is that instead of a setup decision $\Omega_t \in \binary$ 
an integer decision describing the number of containers was used.
Additionally, there is a capacity constraint $\hat{Q}$ for the containers.
The resulting model is:
\LP{ \objd[min]{q,\Omega, \omega}{
		\sum_{t \in T}\left[C^o_{t} \Omega_{t}+\sum_{p \in P} \left[c^o_{tp} \omega_{tp}+c^h_{tp} I_{tp}\right]\right]}}
{   \sto{I_t=I_{t-1,p}+q_{tp}-d_{tp},~ t\in T,p\in P,}
	\sto{\sum_{p\in P} q_{tp}\leq \hat{Q} \Omega_{t},~t\in T}}
A local search heuristic was proposed to solve this problem.
The performance was tested using randomly (uniform, normal) generated demands.
Products varied between 3 and 10 and periods varied between 4 and 18.

\citet{minner2009comparison} introduced three local search heuristics for the multi-product dynamic lot-sizing with limited warehouse capacity and provided a MIP implementation. 
The heuristics are only a few percentage points away from the best solution found.

\begin{figure}[htbp]\centering
	\includegraphics[width=\columnwidth, height=4cm, keepaspectratio]{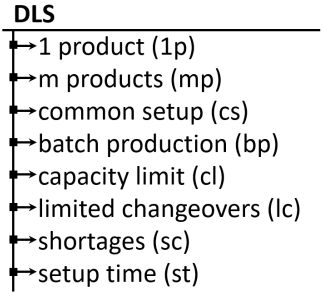}
	\caption{DLS classification}\label{fig:DLS-classification}
\end{figure}

Figure \ref{fig:DLS-classification} shows the proposed classification scheme.
For instance, DLS-mp-cs-cl represents a Dynamic Lot Sizing model with multiple products, common setup and capacity limit.

\subsection{ELSP}
The Economic Lot Sizing Problem (ELSP) is closely related to DLS models.
In general, ELSP attempts to find cycles given constant demand rates.
However, there are several ELSP formulations that allow varying or stochastic demand.
Another main difference is that shortage costs are usually considered within the ELSP.

I have reviewed and classified the ELSP in \cite{Garn2015754}
and the interested reader is referred to this source. 
However, I will add a few reference, which are relevant in regard to MIP formulations.
One of the first and most well-known ELSP formulations is by \citet{bomberger1966dynamic}.
He presented a dynamic programming approach
using the Bellman equation to solve a lot scheduling problem. 
This ELSP does not allow shortages and is similar to the DLS,
but instead of having dynamic demand a constant demand rate per product is given.
In Section \ref{sec:Case-Studies} I will introduce a MIP that offers a solution 
for their classic test-instance.
Another influential paper that took the ELSP into the MIP domain was by
\citet{elmaghraby1978economic}, who discussed the ELSP and gave analytical insights.
The paper used the Bomberger data to demonstrate the concepts.
There are a few more papers, which look at the ELSP using MIP.
\citet{goyal1975scheduling} gave another algorithm to find a solution to the ELSP.
Again, the Bomberger data was used, but credited to \cite{madigan1968scheduling}.
\citet{cooke2004finding} find effective schedules for the economic
lot scheduling problem using a simple mixed integer programming approach.
\cite{salvietti2008profit} look at a profit-maximising economic lot scheduling problem with price optimisation.
A column-generation method was used in combination with a LP to solve the problem.
\cite{sun2010economic} consider the economic lot scheduling problem under extended basic period and power-of-two policy.
Power-of-two means that multiples of twos extend the basic period.

In the above literature review the main focus was on model/problem formulations 
as Mixed Integer Programs.
Now, the focus will shift to how these MIP can be solved in general.

\subsection{Solution methods}
Exact algorithms for the DLS models and ELSPs can be placed into four main categories: 
dynamic programming, branch-and-bound, branch-and-cut and Dantzig-Wolfe decomposition.
\cite{williams1974structured} gave a structured linear programming model in the food industry, where the
Dantzig-Wolfe decomposition was regarded to be ``worthwhile'' for a multiperiod model.

A procedure to tighten bounds was explained in \cite{balas1975facets}. 
The main idea is to find minimum covers for constraints.
$S \subseteq \setn$ is called a \textit{cover} for $\tp{a}\omega \leq \hat{a}$ with $a\in \real^n, \omega \in \binary^n$
if $\sum_{i \in S} a_i > \hat{a}$.
$Q$ is a \textit{minimum cover} if $\sum_{i \in Q} a_i > \hat{a}$ for all $Q \subset S$.
The interesting consequence \cite[Theorem 3]{balas1969structure} is that 
$\sum_{i\in S} \omega_i \leq |S|-1 \equi \tp{a}\omega \leq \hat{a} $.
That means the number of decision variables can be reduced.
This is desirable, because the number of decision variables increases the solution space.
For instance, assume there are $n$ binary variables then there are $2^n$ solutions.
This emphasises the importance of reducing the number of variables.
However, branch-and-bound algorithms cut off many of the $2^n$ solutions,
which are infeasible or too large.

An interesting question is whether and how significant it is to use heuristics
prior to starting exact methods.
Most IP engines use heuristics such as rounding and diving prior to running the branch-and-bound algorithm.
Usually those heuristics  cover the whole solution space.
\cite{linderoth1999computational} review several of the recent IP solution methods and their computational performance.

Recently, an interesting heuristic approach was introduced that focuses on essential
decision variables and then extends the solution space systematically.
For more details see the kernel optimisation methods introduced in 
\cite{guastaroba2014heuristic}, \cite{guastaroba2017adaptive} and \cite{angelelli2010kernel},
although they were used for other problem classes.
This should not be confused with the kernel squeeze feature where 
decision variables with zero value and a zero dual solution cause pivots to reduce the solution space.
This helped \cite{mcbride1988integer} to find solutions for 
the optimal production line selection problem.

Runtime issues with exact algorithms encouraged researchers to propose heuristical approaches
not connected at all to exact approaches.
Greedy-add \citep{federgruen1994joint}, greedy-drop, 
extended-Silver-Meal \citep{silver1973heuristic}, 
general-part-period-balancing \citep{iyogun1991heuristic} and many more are such heuristics.

The above exact solution methods are typically implemented/used within optimisation engines.
Next, I will mention popular Solvers (optimisation engines) 
and sketch out how they work.

\subsection{Optimisation engines \label{sec:OptimisationEngines}}
The models and case studies in this paper were solved using the following three optimisation engines (solvers):
\begin{itemize}
	\item Matlab's optimisation engine,
	\item math.smartana.org optimisation engine: GLPK,
	\item OpenSolver's optimisation engine: COIN-OR CBC linear solver.

\end{itemize}

Matlab's optimisation engine uses the following approach to solve MIPs:
\begin{enumerate}
	\item Preprocess LP by removing redundant variables and constraints;
	\item Solve MIP as relaxed LP;
	\item Preprocess MIP ty tightening bounds, removing redundant constraints, strengthening constraints and fixing integer variables;
	\item Generate cuts such as MI rounding cuts, Gomory cuts, clique cuts, cover cuts, flow cover cuts, strong Chvatal-Gomory cuts, etc.;
	\item Run heuristics to find integer feasible solutions 
	such as rounding, ``diving'' - partial depth-first search \citep{witzig2019conflictdriven}, relaxation induced neighbourhoods \citep{danna2005exploring};
	\item Run branch-and-bound algorithm.
\end{enumerate}

The GLPK (Gnu Linear Programming Kit) performance was analysed in \cite{pryor2011faster}.

The CBC (Coin-or branch and cut) uses LP relaxations, heuristics and cut generations.
Cut generation include rounding cuts, Gomory cuts, clique cuts, Knapsack cover, flow cover and many more.
Hint: For practitioners the OpenSolver Add-In for Excel may be of interest. 
This solver uses the CBC library by default but offers access to other commercial solvers such as Gurobi.

The most popular commercial solvers are Cplex, Xpress and Gurobi.
\cite{meindl2012analysis} analysed the performance of several commercial, free and open-source solvers.
Here, CBC was shown to perform better than GLPK and LP\_Solve.
There are several challenges with measuring performance such as
Solvers being trained for certain test-instances and their variability. 
The performance variability has been discussed in \cite{lodi2013performance}.

This section provided several DLS MIP formulations.
It mentioned methods that can solve these MIP,
and provided a brief insight into optimisation engines used to solve DLS models and ELSP.
The next section introduce a novel general DLS-mp MIP.
This will be followed by case studies demonstrating the performance and limitations.
The Bomberger case study shows how the DLS approach can be used to solve the ELSP.
The conclusion summarises essential findings. 

\section{MIP for the general DLS-mp \label{sec:MIP}}
In the previous section several DLS models were introduced.
Each had a certain specialisation.
The here introduced DLS-mp unifies several of them.
This model introduces several new constraints 
that take into account production limitations during a time period.

There are $m$ products $P=\sett{m}$.
The number of time-slots (weeks) is $T=\sett{n}$.
I will use the index $t$ for time and $p$ to denote a product.
Generally, $t$ describes a period (basic-period, time-window and time-slot are used synonymously).
Without restricting generality, I will use ``week'' to describe $t$, e.g. hours, days and months are equally valid.
There are three types of cost considered: holding, shortage and setup cost.
The holding costs for products are $c^h \in  \real_{+}^{m}$.
The shortage costs are $c^s \in \real_{+}^{m}$, and
changeover costs are $c^o \in  \real_{o}^{m}$. 

Let $I_0 \in  \real^m_+$ be the vector of initial stock, i.e. $I_{0p}$ is the initial stock for product $p$.
The stock for all products during all periods is $I \in \real_{+}^{n \times m}$, 
i.e. $I_{tp}$ represents the stock in period $t$ for product $p$.
The shortage quantity is $S \in \real_{+}^{n \times m}$ and the initial shortage is $S_0 \in \real^m$.
Note that $I_{tp}>0 \ie S_{tp}=0$ and $S_{tp}>0 \ie I_{tp}=0$.
The decisions whether to produce $p$ in $t$ is given by $X=\left(x_{tp}\right) \in \binary^{n \times m }$.
The decisions about the corresponding quantity of $p$ to produce in $t$ are captured in $Q = \left(q_{tp}\right) \in \mathbb{R}_{+}^{n \times m}$. 
Note that $x_{tp}=1 \ie q_{tp}>0$.
The quantities are produced in batches of size $\hat{b}_{p} \in \real_+$ dependent on the product $p$.
Furthermore, the total amount that can be produced during a time-slot is limited to $\hat{q}$.
It should be mentioned here, that these limitations resulted from existing practical requirements.
More flexibility can be introduced by using $\hat{q}_{tp}$ to set a production quantity target for each period and each product.
However, I will refrain from this here.

The lists below summarise the indices, parameters and decision variables explained previously.

\noindent Dimensions and indices:\\
\begin{tabular}{ll}
	$m$  & Number of products \\
	$P$  & Product set \\
	$p$  & Product index \\
	$n$  & Number of time-slots/production windows \\
	$T$  & Time-slot set \\
	$t$  & Time-slot index \\
\end{tabular}%

\noindent Parameters:\\
\begin{tabular}{ll}
	$c^h_p$ & Holding cost for product $p$ (per unit per $t$) \\
	$c^s_p$ & Shortage cost for product $p$ (per unit per $t$) \\
	$c^o_p$ & Changeover cost from any product to product $p$ \\
	$\hat{b}_p$ & Batch size of product $p$ for $t$ \\
	$\hat{q}$ & Production limit during $t$ \\
	$\hat{x}$ & Max. number of products produced during $t$ \\
	$\hat{\omega}$ & Max. number of changeovers during $t$ \\
\end{tabular}%

\noindent Quantities:\\
\begin{tabular}{ll}
	$I_{0p}$ & Initial stock for product $p$ \\
	$I_{tp}$ & Stock/inventory for product $p$ per $t$ \\
	$S_{0p}$ & Initial shortage for product $p$ \\
	$S_{tp}$ & Shortage for product $p$ per $t$ \\
\end{tabular}%

\noindent Decision variables:\\
\begin{tabular}{ll}
	$X_{tp}$ & Decision to produce product $p$ in $t$ \\
	$Q_{tp}$ & Quantity of product $p$ produced in $t$  \\
	$\omega_{tp}$ & Changeover to product $p$ at $t$ \\
\end{tabular}%

\noindent The objective is to minimise the total costs:
\lequ{objwb}{
	\begin{split}
		&\min \text{Cost} =\min \set{ \text{holding} + \text{shortage} + \text{setup} }=\\
		=&\min_{I,S,\omega} \set{\sum_{p\in P} \left[ c^h_p \sum_{t \in T} I_{tp} +  c^s_p \sum_{t\in T} S_{tp} + c^o_p \sum_{t\in T} \omega_{tp} \right] }
	\end{split}
}

The inventory needs to fulfil the balance equation:
\lequ{inventoryBalance}{\left(I_{t-1,p}-S_{t-1,p}\right)+q_{t p}-d_{t p} \stackrel{!}{=} I_{t p}-S_{t p},}
where $t \in T$ and $p \in P$. Note that $I_{0p}$ and $S_{0p}$ are known.
The previous period's stock/shortage is topped up by the production quantity in the current period
and reduced by the current demand. This sets the period's stock/shortage.
For instance, assume last week's stock was 10 units. 
This week's production quantity is 40 units and there is a demand of 20.
Hence, the current week's stock is $10+40-20 = 30$.

In the ELSP only one product is permitted to be produced per period (e.g. week):
\lequ{nbProductsProducedOne}{\sum_{p \in P} x_{tp} \meq 1, ~t \in T.}
However, there are other interpretations: 
(1) although there is only one production line several products can be produced during a sufficiently large period;
(2) there are several production lines available to produce a single batch during a period;
(3) several production lines and several products are produced simultaneously during a period.
To take care of the above interpretations constraint (\ref{eq:nbProductsProducedOne}) is
changed to 
\lequ{nbProductsProducedMany}{\sum_{p \in P} x_{tp} \leq \hat{x}, ~t \in T,}
where $\hat{x}$ is the maximum number of different products produced simultaneously.
Furthermore, the production limit constraint (\ref{eq:productionLimit}) ``kicks in''.

The amount produced is exactly the same as the batch size:
\lequ{batchLimit}{q_{tp} \meq  x_{tp} \hat{b}_{p}, ~t\in T, p \in P.}
Here, it can make sense to allow the production of multiple batch sizes,
which would require the decision variable $x_{tp}$ to be in $\setN_0$ rather than $\binary$.

The total amount produced cannot exceed the time-slot's (e.g. weekly) production limit:
\lequ{productionLimit}{\sum_{p \in P} q_{tp} \leqslant \hat{q},~t \in T.}

Now, the ``logical'' IP constraints that identify changeovers are given by:
\lequ{changeoverConstraints}{
	\begin{aligned}
		&x_{t+1,p}-x_{tp} \geq 2\omega_{t+1, p}-1,\\
		&x_{t+1,p}-x_{tp} \leq 2\omega_{t+1, p},\\
		& t \in \sett{n-1},~ p \in P.
	\end{aligned}
}
The initial setups (changeovers) are identical to whether the product is produced:
\lequ{initialSetupConstraints}{
	\omega_{1, p} \meq x_{1,p} ,~ p \in P.
}

\begin{figure}[htbp]\centering
	\includegraphics[width=\linewidth]{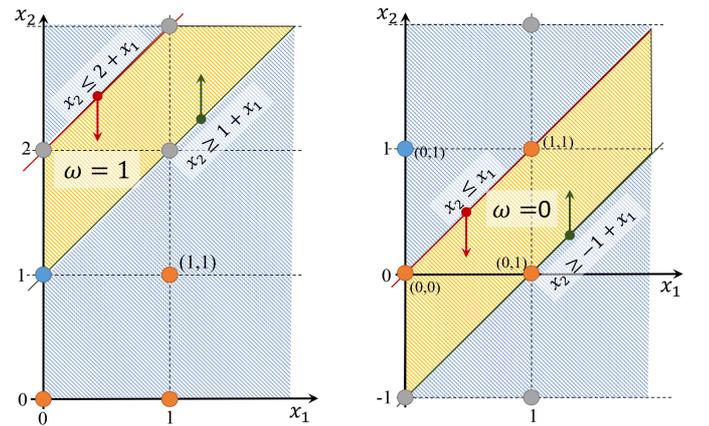} 
	\caption{Logical constraints}\label{fig:logicalconstraints}
\end{figure}

The logical constraints require validity checks. 
The problem is simplified by considering one transition.
Let $x_1$ and $x_2$ be the current and next production state. 
If $x_1=1$ then something is produced otherwise not. 
A setup-changeover happens when $(x_1,x_2)=(0,1)$. This is reflected with $\omega = 1$. 
The discontinuation of producing $(x_1,x_2)=(1,0)$ is also a changeover.
However, in this formulation it is assumed that the discontinuation-changeover activities and associated costs 
are integrated within the setup-changeover, i.e. $\omega = 0$.
If $(x_1,x_2) \in \set{(0,0),(1,1)}$, then no changeover takes place, i.e. $\omega = 0$.
This represents the logical operation $\neg x_1 \land x_2$ or $x_1<x_2$ shown in Table \ref{tab:lop}.
\qtable{Logical operation $\neg x_1 \land x_2 \ie \omega$.}{lop}{r|r|r|r}
{$x_1$ & $x_2$ & $\neg x_1$ & $\omega$}
{0     & 0   & 1  & 0 \\
	0     & 1   & 1  & 1 \\
	1     & 0   & 0  & 0 \\
	1     & 1   & 0  & 0 \\}
The corresponding logical constraint is:
\qbox{\tcb{If} $x_1<x_2$ \tcb{then} $\omega = 1$ \tcb{else} $\omega=0$.}
This can be written as IP constraints (motivated by Figure \ref{fig:logicalconstraints}):
\lequ{logicalConstraintsSimple}{
\begin{aligned}
&x_2-x_1 \geq 2\omega-1,\\
&x_2-x_1 \leq 2\omega.
\end{aligned}}
In order to show that this is a valid formulation 
the logical constraints must be equivalent to the IP-constraints.

Suppose that $(x_1, x_2,\omega)$ is a feasible solution for the IP formulation.
If $\omega=1$ then $x_2-x_1 \geq 1$ and $x_2-x_1 \leq 2$.
The left of Figure \ref{fig:logicalconstraints} shows the constraints and  feasible region.
This implies that $(x_1,x_2) \meq (0,1)$, because this is the only feasible solution fulfilling both constraints. 
The other points $\set{(0,0),(1,0),(1,1)}$ are infeasible solutions.
If $\omega=0$ then $x_2-x_1 \geq -1$ and $x_2-x_1 \leq 0$.
The right of Figure \ref{fig:logicalconstraints} shows the constraints and feasible region.
This implies that $(x_1,x_2) \in \set{(0,0),(1,0),(1,1)}$, because a feasible solution was assumed. 
The other point $(x_1,x_2)=(0,1)$ is an infeasible solution.
This proves the validity of (\ref{eq:logicalConstraintsSimple}) 
and thus the validity of (\ref{eq:changeoverConstraints}).

$\omega_{tp}$ represents a changeover.
If only one product can be produced during a period, then 
constraint (\ref{eq:nbProductsProducedOne}) will suffice, 
and ensure only one changeover.
However, if multiple products can be produced simultaneously,
then constraint (\ref{eq:nbProductsProducedMany}) is used.
However, it may be necessary to limit the number of changeovers
to $\hat{\omega} = (\hat{\omega}_t) \in \real_+^n$:
\lequ{changeOverLimit}{\sum_{p\in P} \omega_{tp} \leq \hat{\omega}_t,~t\in T.}
This can happen if the number of resources available to 
execute the changeovers is constrained.






\section{Case studies \label{sec:Case-Studies}}

Two cases studies are presented.
The first demonstrates the functionality of the general DLS model introduced in the previous section as DLS-mp-bp-cl-lc-sc system (see Figure \ref{fig:DLS-classification}). 
The second one illustrates how a classic ELSP 
can be solved using a DLS-mp-qc-st system approach.

\subsection{Food manufacturer}
The case study, which was introduced in \cite{Garn2015754} will be used here.
For the convenience of the reader a brief summary is given.
A company produces several SKUs (stock keeping units) and requires a production schedule.
The case study consists out of $m=10$ SKUs produced over a period of two years.
The demand data is given on a weekly basis ($n=104$), i.e. $D\in \real^{104\times 10}$.
The demand was 404k packages in the first year and 608k in the second packages.
Figure \ref{fig:total_demand_per_week} shows the total weekly demand.
\begin{figure}[htbp]\centering
	\includegraphics[width=\linewidth]{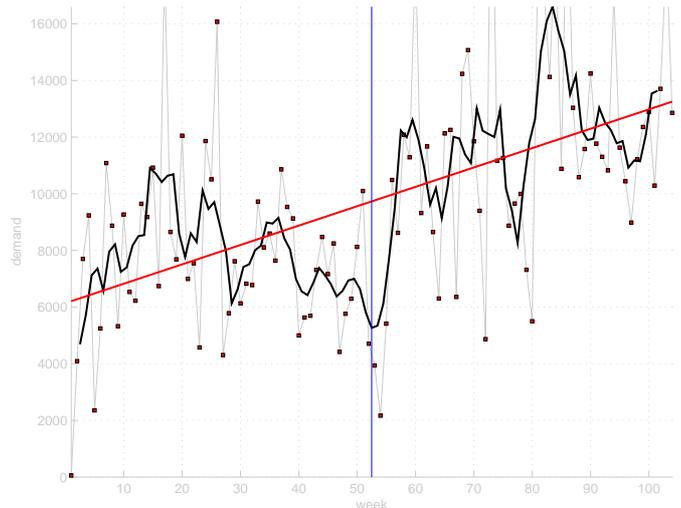} 
	\caption{Total demand per week.}\label{fig:total_demand_per_week}
\end{figure}
This demand is assembled via 10 individual SKUs shown in Figure \ref{fig:total_demand_per_week}.
\begin{figure*}[htbp]\centering
	\includegraphics[width=\linewidth]{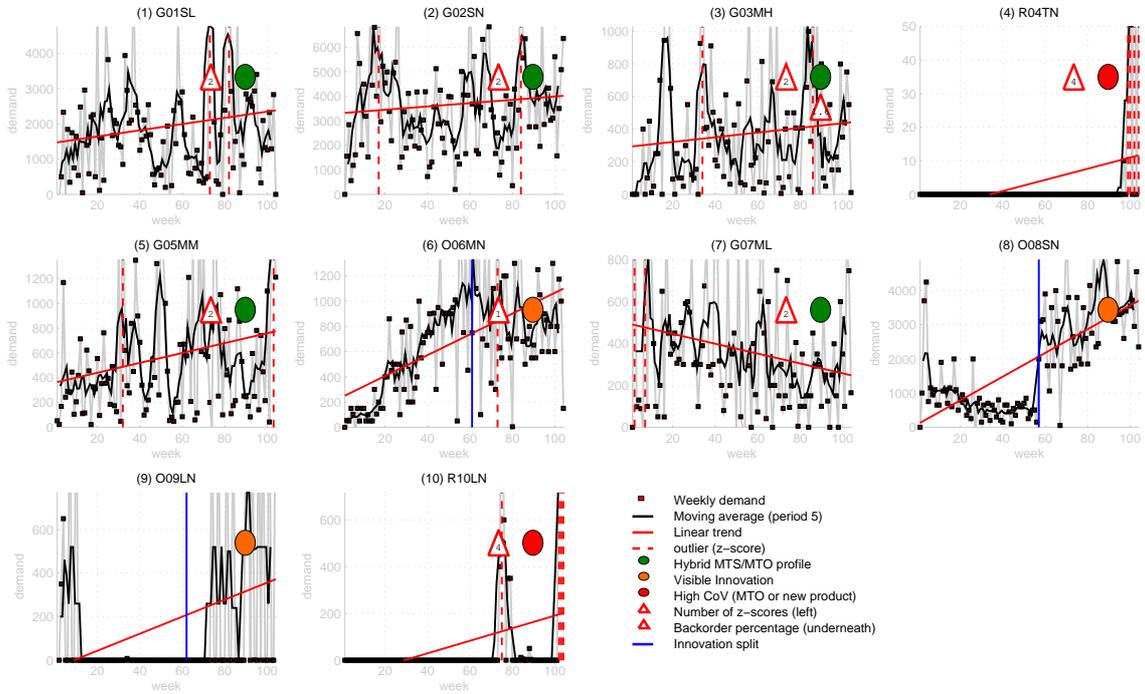} 
	\caption{Demand and characteristics of each SKU.}\label{fig:total_demand_per_week}
\end{figure*}

Holding costs varied between \$1.10 and \$3.00 per package per week.
The shortage costs are in the range \$4.60 and \$8.60 per package per week.
All packages can be back-ordered.
The  unit-costs for all products are fixed over the entire two years.
Changeover (setup) costs are between \$5,200 and \$8,300.
They only depend on the product that is produced, i.e. independent of which product was produced before.
The detailed costs, batch-sizes and aggregated demands are shown in Table \ref{tab:cost-factors}.
\input{tab-cost-factors}
The table illustrates that three products (G02SN, G01SL and O08SN) 
have a demand of 77.2\%.
Hence, in case of complexity issues the focus should be on these three products.
Accordingly, the products were ordered by their total demand.
The weekly production capacity of 5,000 units is the same as the maximum batch-size.
Furthermore, all other batch sizes suggest that producing more than one product a week makes sense.

\begin{figure}[htbp]\centering
	\includegraphics[width=\linewidth]{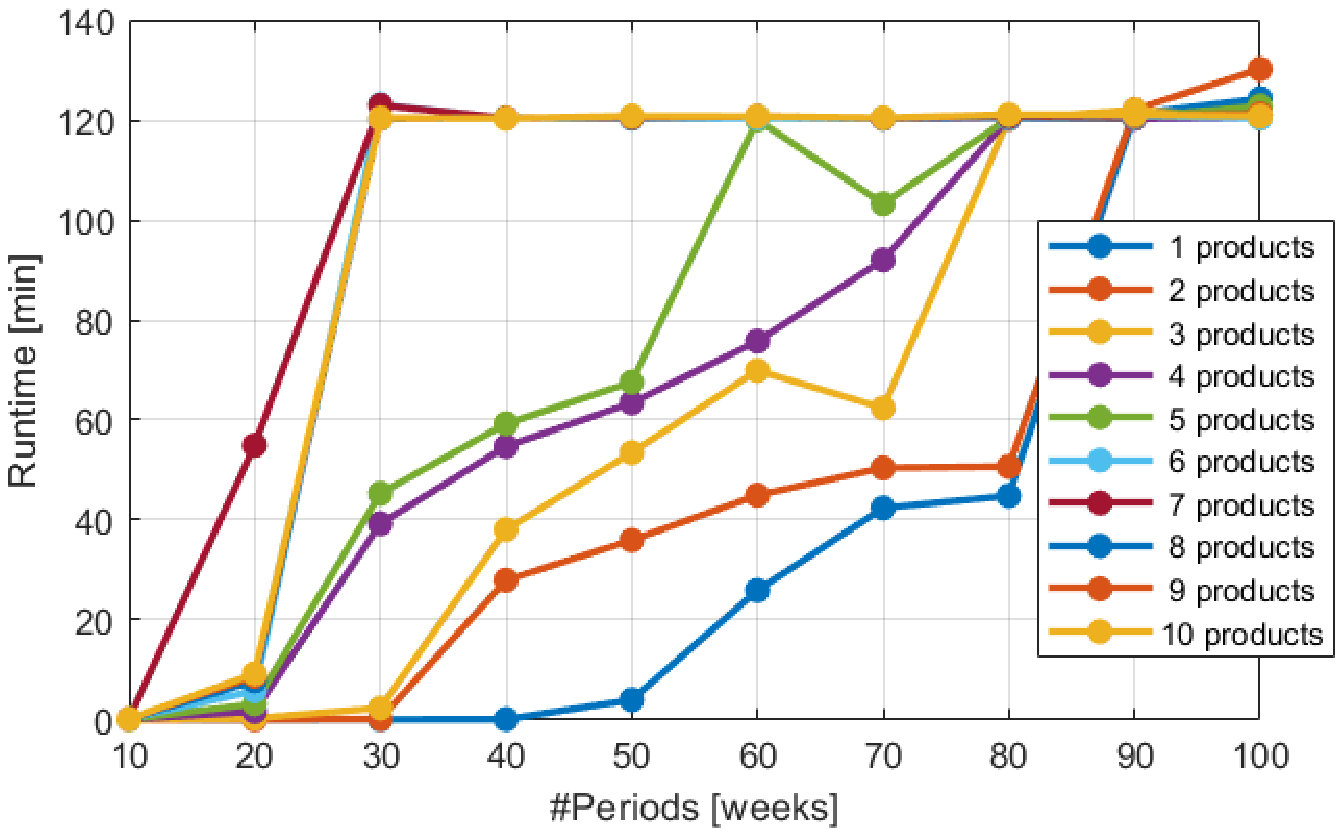} 
	\caption{Runtime for test-instances of case study.}\label{fig:cs-runtime}
\end{figure}

The case study's complexity was increased by 10 weeks and one product. 
Hence, 100 test-instances based on this one case study were created.
The Matlab optimisation engine (see Section \ref{sec:OptimisationEngines}) was used 
to solve the MIP explained in Section \ref{sec:MIP}.
The time to solve the MIP problem was limited to two hours.
46 out of 100 test-instances obtained an optimal solution within the set time limit.
Figure \ref{fig:cs-runtime} shows the total runtime, i.e. time to formulate problem plus time for solving it. 
For instance, the runtime to find the production schedule for the first 50 weeks 
and the top three products is 53.5 minutes.
This involves the usage of a virtual machine with 4 cores with 2GHz allocated Intel CPUs 
and 1GB RAM (out of 8GB) available.
The total number of decision variables used was $|T|\times|P|\times 4 = 50\times 3\times 4$,
where the last 4 represents the number of decision variable matrices ($\set{Q,X,\omega,I,S}$).
As expected the MIP finds optimal solutions with moderate complexity 
within the two-hours time window.
It is interesting to observe that the optimisation engine fails faster to find optimal solutions,
when the number of products increases rather than when time increases.
For instance, for up to 4 products and 70 weeks ($280\times 4$ variables) an optimal solution is still found.
However, at 6 products and 30 weeks ($180\times 4$ variables) only a sub-optimal solution is found.

\begin{figure}[htbp]\centering
	\includegraphics[width=\linewidth]{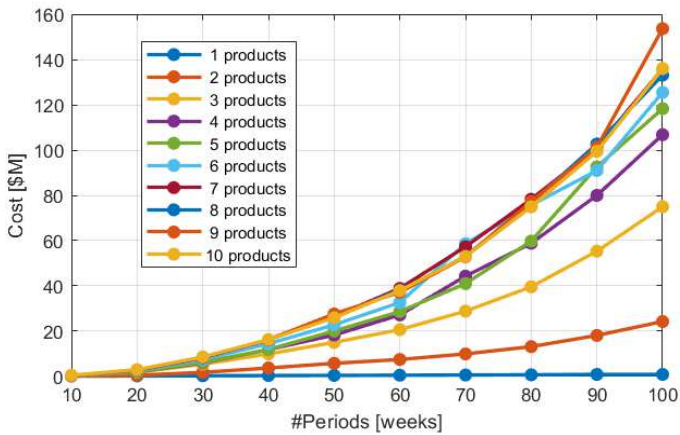} 
	\caption{Cost for test-instances of case study.}\label{fig:cs-cost}
\end{figure}

Figure \ref{fig:cs-cost} shows the cost increase over the weeks and products produced.
If the demand is roughly uniform, one would expect a linear growth of costs.
However, the demand increases as shown in Figure \ref{fig:total_demand_per_week}.
Moreover, there are outliers, which put additional strain/costs on the system 
- see \cite{Garn2015} for details.
The above explains why cost are increasing apparently exponential.

Next the effect on runtime and solution quality will be considered if an initial solution is provided.
Here, the initial solution for $m=1$ is not to produce anything, 
and for $m>1$ to use the solution found for $m-1$ and not to produce product $m$.
This is a similar idea as introduced in \cite{Garn2015754} in the iterative schedule generation algorithm.
The following results do not support the results expected intuitively.
Runtime savings would have been expected. 
Especially, since the time to find an initial solution was not taken into account.
Surprisingly,  optimal solutions for only 28\% of the test-instances were found, 
when an initial solution was provided. Previously 46\% optimal solutions were determined.
62\% of the instances had longer runtimes.
Solution quality is 51\% of the time worse than previously.
Since the experiments were run on shared VM, 
the unexpected results may be due to shared memory and CPUs being affected,
or because of an implementation issue within Matlab's optimisation engine.

Overall, the case study demonstrated that optimal solutions 
can be found for most test instances. 
The number of products influences the runtime significantly.
Since the top three products dominate the economic benefits, this issue can be mitigated.

\subsection{Bomberger}
The classic test data from \cite{bomberger1966dynamic} was used.
This ELSP focuses on production cycles.
Daily demand for each product per day $d_p$ is given and assumed to be constant (over the entire time-horizon).
The test-instance investigated periods were 20, 35 and 40 days.
Practically, it does not make sense to plan the schedule for more than two years.
A year has $y=240$ days in this problem instance.
That means, if $\tau=20$ then there are $n=24$ periods to consider.
The period $\tau$ are long in comparison to the setup time $\tau^o$, 
which would allow us to neglect the setup times.
However, it is equally convenient to consider the setup time 
within the maximum quantity producible during $\tau$.
It appears reasonable to require setup costs in each period,
even when there is a continuation of the production into the next period.
This is due to the extensive length of the production period.
The production rate is given in units per day.
Hence, the maximum is $\hat{q}_p = (\tau-\tau^o) b_p$, 
where $b_p$ represents the production rate for product $p$ (Note: ensure identical units such as days for $\tau$ and $\tau^o$).
Assume that there is a minimum amount that has to be produced, e.g. $\check{q}_p = \frac{b_p}{8}$ (i.e. 1 hour runtime).
The amount produced is captured in the decision variable $q_{tp} \in \real_{+}$.
The number of days a product is produced is abbreviated by $x_{tp}\in \real_{+}$, i.e. fractional days are allowed.
The fact whether a product is produced is equivalent to a changeover/setup and represented by $\omega_{tp}\in \binary$.
The daily holding cost is defined as $f=0.10$  multiplied by the cost of
labour and material (i.e. piece cost $c^m$). 
Hence, the holding cost per unit during $\tau$ is $c^h = 0.1 c^m \tau$. 
This interpretation introduces imprecision, 
e.g. if something is produced on the last day of $\tau$ then the entire period needs to be paid.
Yet, in real-world scenarios it is common to arrange a minimum payable period (e.g. one week) 
and a minimum quantity (e.g.one pallet) as payment to third-party warehouses.
Usually initial inventories $I_{0p} \in \real_{+}$ are present. 
However, in this test-instance $I_{0p}=0$.
Furthermore, the Bomberger test-instance does not permit shortages.
The above allows us to give a simplified formulation 
of the general DLS introduced in Section \ref{sec:MIP}:
\LP{\label{MIP:Bomberger} 
	\objd[min]{I,\omega}{\sum_{p \in P}\sum_{t \in T} (c^h_p I_{tp} + c^o_p \omega_{tp})}}
{
	\sto{I_{tp} = I_{t-1,p} + q_{tp} - d_{p} \mcomment{inv. balance}}
	\sto{\sum_{p\in P} (x_{tp} + \omega_{tp}\tau^o_p) \leq \tau \mcomment{max time}}
	\sto{ x_{tp} b_{p} = q_{tp} \mcomment{production rate}}
	\sto{q_{tp}\leq \hat{q}_p \mcomment{quantity limit}}
	\sto{q_{tp} \leq \hat{q}_p \omega_{tp} \mcomment{if $q_{tp}$ then $\omega_{tp}$}}
	\sto{q_{tp} \geq \check{q}_p \omega_{tp} \mcomment{if $\omega_{tp}$ then $q_{tp}$}}	
	\sto{I_{tp}, q_{tp}, x_{tp}\in \real^+,~\omega_{tp} \in \binary \mcomment{variables}.}
}
It is assumed that $t \in T$ and $p\in P$, 
where $T=\sett{n}$ is the number of periods during the schedule horizon ($n=\frac{y}{\tau}$);
and $P=\sett{m}$ represents $m$ products.
All factors have to be scaled to have consistent units.

Upper bounds for the optimisation variables are $I_{tp} \leq n \times \max_p \set{d_p}$,
$q_{tp} \leq \max_p \set{\hat{q}_p}$ and $x_{tp}\leq \tau$.  

The parameters: setup cost $c^o$, piece cost $c^m$, production rate $b$, demand rate $d$ and setup time $\tau^o$ 
are shown in Table \ref{tab:BombergerData}.
\begin{table}[htbp]
	\centering
	\caption{Bomberger data.}
	\adjustbox{width=\columnwidth}{ \begin{tabular}{c|cc|cc|c}
		\toprule
		Part Nb & Setup Cost & Piece Cost & Prod. Rate & Demand & Setup Time \\
		& [\$] & [\$/unit] & [units/day] & [per day] & [hours] \\
		\midrule
		1    & 15   &        0.0065  &          30,000  &          400  & 1 \\
		2    & 20   &        0.1775  &           8,000  &          400  & 1 \\
		3    & 30   &        0.1275  &           9,500  &          800  & 2 \\
		4    & 10   &        0.1000  &           7,500  &       1,600  & 1 \\
		5    & 110  &        2.7850  &           2,000  &            80  & 4 \\
		6    & 50   &        0.2675  &           6,000  &            80  & 2 \\
		7    & 310  &        1.5000  &           2,400  &            24  & 8 \\
		8    & 130  &        5.9000  &           1,300  &          340  & 4 \\
		9    & 200  &        0.9000  &           2,000  &          340  & 6 \\
		10   & 5    &        0.0400  &        15,000  &          400  & 1 \\
		\bottomrule
	\end{tabular}}%
	\label{tab:BombergerData}%
\end{table}%

In general, the assertion:$\frac{d_p}{b_p} \leq 1$ must hold for a feasible solution.
Product (part nb, item) 8 has a $\frac{d_p}{b_p} = \frac{340}{1300}\approx 0.262$ fraction, which is the maximum.
This value can be interpreted that it takes a quarter of a day to produce the required demand for that day.
If there exists a feasible solution, then
\lequ{BB-necessary-condition}{\sum_{p\in P} \left(\frac{d_p\tau-I_0}{b_p} \right) \leq \tau,}
because the demand must be fulfilled already in the first cycle.
In the Bomberger example the sum is 21.4 days,
which is well below the $\tau$ (35 and 40 days) used in Bomberger's paper.
For $\tau=20$ we will assume an initial inventory of $I_0=\frac{d\tau}{2}$.
Bomberger stated a similar necessary condition for repetitive schedules 
but allowing different periods for each product.

Solving this MIP (\ref{MIP:Bomberger}) gives a daily cost of \$444 for $\tau=20$.
This is a higher solution value than Bomberger's. 
This may be due to the used period length.
For instance, a period of $\tau=80$ and $\tau=120$  finds a solution value which is \$130 and \$7.33 (with no initial inventory).
Here, the holding and setup cost drop significantly. 
Only as much is produced as there is consumed. 
Bomberger's cycle lengths are up to 400 days, 
which appear to be unreasonably high for a manufacturing environment.
Furthermore, demand was assumed as uniform, which will not happen in real world scenarios.

A period of five and ten working-days is examined next. 
An initial inventory for five days is provided for all products.
The minimum quantity produced is set to a two hours runtime, i.e. $\check{q} =\frac{1}{4}b$.
The average daily costs are \$241.78 (relative gap: 36.8\%, after 90min runtime) for $\tau=5$ 
and \$109.88 (optimal) for $\tau=10$.
$\omega_p$ reveals almost cyclic solutions.
The average number of days between setups is:
\smat{38.8& 5.0& 5.0& 5.0& 9.8& 21.7& 31.8& 5.0& 5.0& 9.8} for $\tau=5$
and \smat{20.0& 10.0& 10.0& 10.0& 10.0& 20.4& 30.0& 10.0& 10.0& 10.2} for $\tau=10$.

This case study demonstrated that the ELSP can be solved using a typical DLS modelling approach.

\section{Conclusions}
DLS has received a fair amount of attention in the past,
but several essential constraints were omitted
such as production limits and shortages.
Sometimes, unnecessary (from an industry point of view) constraint specifications  were introduced.
For instance, holding costs varying at every time period.
Interesting first steps into conjoint sales price considerations have appeared,
but more has to be done for it to be useable in practice.
This work has unified several of the previous studies and suggested alternative formulations and approaches when appropriate.

The main contributions of this work are the classification of DLS models 
and a novel DLS MIP formulation,
which is suitable for agile manufacturing in an industrial environment
allowing shortages and introducing production \& changeover limits 
(i.e. a DLS-mp-bp-cl-lc-sc system).
In the context of the MIP the validity of changeover constraints was proven.
The capacitated DLS models for multiple products and ELSPs have been implemented as MIPs using Matlab's optimisation engine.
A case-study consisting out of 200 test-instances was executed.
Interesting run-time and cost insights were gained, such as product variety causes more performance issues than longer time periods.
The classic ELSP Bomberger test instance was solved using a MIP.

To conclude, the here introduced MIPs can be used for a variety of industry solutions in the area of agile manufacturing.
A unified approach in formulating and classifying DLS models and ELSP was proposed.

\section{Appendix}
Table \ref{tab:cs-runtime} and \ref{tab:cs-cost} show the runtimes and costs for test-instances resulting from the food manufacturer case study respectively.
\input{tab-cs-runtime}
\input{tab-cs-cost}

\bibliographystyle{elsarticle-harv} 
\bibliography{literature}
\end{document}

%% file: tab-cost-factors.tex
\begin{table*}[htbp]
	\centering
	\caption{Cost factors, batch-sizes and demand.}
	\adjustbox{width=\textwidth}{
	\begin{tabular}{l|cccccccccc}
		\toprule
		Cost & G01SL & G02SN & G03MH & R04TN & G05MM & O06MN & G07ML & O08SN & O09LN & R10LN \\
		\midrule
		Holding &                 1.30  &                 2.40  &              1.50  &            2.20  &              1.10  &              1.60  &              2.30  &                 3.00  &              2.30  &            2.20  \\
		Setup &               5,600  &               7,900  &            5,600  &          6,500  &            7,000  &            7,100  &            5,600  &               7,700  &            5,200  &          8,300  \\
		Shortage &                 7.20  &                 7.60  &              6.60  &            8.10  &              4.60  &              8.60  &              7.30  &                 8.10  &              6.80  &            8.30  \\
		\midrule
		Batch-size &               4,000  &               5,000  &                700  &              100  &            1,100  &            1,300  &                700  &               3,800  &                300  &              100  \\
		\midrule
		Demand &          200,626  &          381,381  &          38,241  &              320  &          59,228  &          70,167  &          38,429  &          199,239  &          17,759  &          6,576  \\
		Demand [\%] &                 19.8  &                 37.7  &                 3.8  &                 -    &                 5.9  &                 6.9  &                 3.8  &                 19.7  &                 1.8  &               0.6  \\
		\bottomrule
	\end{tabular}}
	\label{tab:cost-factors}%
\end{table*}%

%% file: tab-cs-runtime.tex
  \begin{table*}[htbp]
    \centering
    \caption{Case study runtime.}
      \begin{tabular}{c|cccccccccc}
      \toprule
      weeks & P1   & P2   & P3   & P4   & P5   & P6   & P7   & P8   & P9   & P10 \\
      \midrule
      10   &             0.3  &             0.5  &             1.0  &             2.3  &             2.0  &           11.0  &             9.4  &             6.7  &             7.6  &             8.3  \\
      20   &             0.6  &             1.0  &           16.4  &           98.9  &         189.0  &         346.7  &     3,296.0  &         465.3  &         521.7  &         552.2  \\
      30   &             1.0  &             7.3  &         138.8  &         2,354  &         2,718  &         7,388  &         7,380  &         7,219  &         7,222  &         7,228  \\
      40   &             3.2  &         1,675  &         2,284  &         3,282  &         3,555  &         7,226  &         7,226  &         7,231  &         7,228  &         7,223  \\
      50   &         237.0  &         2,156  &         3,207  &         3,805  &         4,054  &         7,224  &         7,232  &         7,228  &         7,235  &         7,251  \\
      60   &         1,552  &         2,695  &         4,194  &         4,547  &         7,206  &         7,230  &         7,249  &         7,249  &         7,249  &         7,247  \\
      70   &         2,543  &         3,025  &         3,744  &         5,522  &         6,198  &         7,221  &         7,222  &         7,224  &         7,227  &         7,230  \\
      80   &         2,689  &         3,040  &         7,215  &         7,220  &         7,226  &         7,229  &         7,242  &         7,254  &         7,262  &         7,268  \\
      90   &         7,278  &         7,329  &         7,325  &         7,219  &         7,237  &         7,247  &         7,256  &         7,261  &         7,261  &         7,266  \\
      100  &         7,461  &         7,816  &         7,225  &         7,232  &         7,373  &         7,228  &         7,267  &         7,266  &         7,275  &         7,249  \\
      \bottomrule
      \end{tabular}%
    \label{tab:cs-runtime}%
  \end{table*}%

%% file: tab-cs-cost.tex
  \begin{table*}[htbp]
    \centering
    \caption{Case study cost.}
      \adjustbox{width=\textwidth}{
      \begin{tabular}{c|cccccccccc}
      \toprule
      weeks & P1   & P2   & P3   & P4   & P5   & P6   & P7   & P8   & P9   & P10 \\
      \midrule
      10   &           59,001  &                 95,016  &               236,662  &                 266,719  &                 282,196  &                 391,420  &                 412,574  &                 469,644  &                 469,644  &                 469,644  \\
      20   &         165,915  &               451,359  &           1,788,770  &              1,943,398  &              2,006,730  &              2,527,422  &              2,788,540  &              2,966,089  &              2,966,089  &              2,966,089  \\
      30   &         249,965  & 1.74M & 5.33M & 5.92M & 6.14M & 7.46M & 8.25M & 8.53M & 8.53M & 8.53M \\
      40   &         322,581  & 3.67M & 9.89M & 11.86M & 11.90M & 14.39M & 16.07M & 16.22M & 16.24M & 16.21M \\
      50   &         399,939  & 5.74M & 15.00M & 18.24M & 20.02M & 22.77M & 26.49M & 26.83M & 27.57M & 25.86M \\
      60   &         480,938  & 7.43M & 20.63M & 27.12M & 28.63M & 32.61M & 38.87M & 37.99M & 37.18M & 37.90M \\
      70   &         577,721  & 9.90M & 28.80M & 44.28M & 40.99M & 58.65M & 57.30M & 53.05M & 52.82M & 53.00M \\
      80   &         657,683  & 13.12M & 39.56M & 58.82M & 59.79M & 75.69M & 78.29M & 76.79M & 77.24M & 74.95M \\
      90   &         842,058  & 18.07M & 55.30M & 79.99M & 92.57M & 90.98M & 101.12M & 102.69M & 101.27M & 99.48M \\
      100  &         844,760  & 24.23M & 74.97M & 106.82M & 118.34M & 125.49M & 135.98M & 133.20M & 153.70M & 136.03M \\
      \bottomrule
      \end{tabular}
    }
    \label{tab:cs-cost}%
  \end{table*}%

%% file: DLS-MIP.bbl
\begin{thebibliography}{36}
\expandafter\ifx\csname natexlab\endcsname\relax\def\natexlab#1{#1}\fi
\providecommand{\url}[1]{\texttt{#1}}
\providecommand{\href}[2]{#2}
\providecommand{\path}[1]{#1}
\providecommand{\DOIprefix}{doi:}
\providecommand{\ArXivprefix}{arXiv:}
\providecommand{\URLprefix}{URL: }
\providecommand{\Pubmedprefix}{pmid:}
\providecommand{\doi}[1]{\href{http://dx.doi.org/#1}{\path{#1}}}
\providecommand{\Pubmed}[1]{\href{pmid:#1}{\path{#1}}}
\providecommand{\bibinfo}[2]{#2}
\ifx\xfnm\relax \def\xfnm[#1]{\unskip,\space#1}\fi
\bibitem[{Absi and Kedad-Sidhoum(2009)}]{absi2009multi}
\bibinfo{author}{Absi, N.}, \bibinfo{author}{Kedad-Sidhoum, S.},
  \bibinfo{year}{2009}.
\newblock \bibinfo{title}{The multi-item capacitated lot-sizing problem with
  safety stocks and demand shortage costs}.
\newblock \bibinfo{journal}{Computers \& Operations Research}
  \bibinfo{volume}{36}, \bibinfo{pages}{2926--2936}.
\bibitem[{Angelelli et~al.(2010)Angelelli, Mansini and
  Speranza}]{angelelli2010kernel}
\bibinfo{author}{Angelelli, E.}, \bibinfo{author}{Mansini, R.},
  \bibinfo{author}{Speranza, M.G.}, \bibinfo{year}{2010}.
\newblock \bibinfo{title}{Kernel search: A general heuristic for the
  multi-dimensional knapsack problem}.
\newblock \bibinfo{journal}{Computers \& Operations Research}
  \bibinfo{volume}{37}, \bibinfo{pages}{2017--2026}.
\bibitem[{Balas(1975)}]{balas1975facets}
\bibinfo{author}{Balas, E.}, \bibinfo{year}{1975}.
\newblock \bibinfo{title}{Facets of the knapsack polytope}.
\newblock \bibinfo{journal}{Mathematical programming} \bibinfo{volume}{8},
  \bibinfo{pages}{146--164}.
\bibitem[{Balas and Jeroslow(1969)}]{balas1969structure}
\bibinfo{author}{Balas, E.}, \bibinfo{author}{Jeroslow, R.},
  \bibinfo{year}{1969}.
\newblock \bibinfo{title}{On the structure of the unit hypercube}.
\newblock \bibinfo{journal}{Management Sci. Research Rep.}
  \bibinfo{volume}{198}.
\bibitem[{Boctor et~al.(2004)Boctor, Laporte and Renaud}]{boctor2004models}
\bibinfo{author}{Boctor, F.F.}, \bibinfo{author}{Laporte, G.},
  \bibinfo{author}{Renaud, J.}, \bibinfo{year}{2004}.
\newblock \bibinfo{title}{Models and algorithms for the dynamic-demand joint
  replenishment problem}.
\newblock \bibinfo{journal}{International Journal of Production Research}
  \bibinfo{volume}{42}, \bibinfo{pages}{2667--2678}.
\bibitem[{Bomberger(1966)}]{bomberger1966dynamic}
\bibinfo{author}{Bomberger, E.E.}, \bibinfo{year}{1966}.
\newblock \bibinfo{title}{A dynamic programming approach to a lot size
  scheduling problem}.
\newblock \bibinfo{journal}{Management science} \bibinfo{volume}{12},
  \bibinfo{pages}{778--784}.
\bibitem[{Br{\"u}ggemann and Jahnke(2000)}]{bruggemann2000discrete}
\bibinfo{author}{Br{\"u}ggemann, W.}, \bibinfo{author}{Jahnke, H.},
  \bibinfo{year}{2000}.
\newblock \bibinfo{title}{The discrete lot-sizing and scheduling problem:
  Complexity and modification for batch availability}.
\newblock \bibinfo{journal}{European Journal of Operational Research}
  \bibinfo{volume}{124}, \bibinfo{pages}{511--528}.
\bibitem[{Cooke et~al.(2004)Cooke, Rohleder and Silver}]{cooke2004finding}
\bibinfo{author}{Cooke, D.L.}, \bibinfo{author}{Rohleder, T.R.},
  \bibinfo{author}{Silver, E.A.}, \bibinfo{year}{2004}.
\newblock \bibinfo{title}{Finding effective schedules for the economic lot
  scheduling problem: A simple mixed integer programming approach}.
\newblock \bibinfo{journal}{International Journal of Production Research}
  \bibinfo{volume}{42}, \bibinfo{pages}{21--36}.
\bibitem[{Danna et~al.(2005)Danna, Rothberg and Le~Pape}]{danna2005exploring}
\bibinfo{author}{Danna, E.}, \bibinfo{author}{Rothberg, E.},
  \bibinfo{author}{Le~Pape, C.}, \bibinfo{year}{2005}.
\newblock \bibinfo{title}{Exploring relaxation induced neighborhoods to improve
  mip solutions}.
\newblock \bibinfo{journal}{Mathematical Programming} \bibinfo{volume}{102},
  \bibinfo{pages}{71--90}.
\bibitem[{Elmaghraby(1978)}]{elmaghraby1978economic}
\bibinfo{author}{Elmaghraby, S.E.}, \bibinfo{year}{1978}.
\newblock \bibinfo{title}{The economic lot scheduling problem (elsp): review
  and extensions}.
\newblock \bibinfo{journal}{Management Science} \bibinfo{volume}{24},
  \bibinfo{pages}{587--598}.
\bibitem[{Federgruen and Tzur(1994)}]{federgruen1994joint}
\bibinfo{author}{Federgruen, A.}, \bibinfo{author}{Tzur, M.},
  \bibinfo{year}{1994}.
\newblock \bibinfo{title}{The joint replenishment problem with time-varying
  costs and demands: Efficient, asymptotic and $\varepsilon$-optimal
  solutions}.
\newblock \bibinfo{journal}{Operations Research} \bibinfo{volume}{42},
  \bibinfo{pages}{1067--1086}.
\bibitem[{Garn(2018)}]{Garn2018Dec}
\bibinfo{author}{Garn, W.}, \bibinfo{year}{2018}.
\newblock \bibinfo{title}{{Introduction to Management Science: Modelling,
  Optimisation and Probability}}.
\newblock \bibinfo{publisher}{Smartana Ltd}.
\bibitem[{Garn and Aitken(2015a)}]{Garn2015754}
\bibinfo{author}{Garn, W.}, \bibinfo{author}{Aitken, J.},
  \bibinfo{year}{2015}a.
\newblock \bibinfo{title}{Agile factorial production for a single manufacturing
  line with multiple products}.
\newblock \bibinfo{journal}{European Journal of Operational Research}
  \bibinfo{volume}{245}, \bibinfo{pages}{754 -- 766}.
\bibitem[{Garn and Aitken(2015b)}]{Garn2015}
\bibinfo{author}{Garn, W.}, \bibinfo{author}{Aitken, J.},
  \bibinfo{year}{2015}b.
\newblock \bibinfo{title}{{Splitting hybrid Make-To-Order and Make-To-Stock
  demand profiles}}.
\newblock \bibinfo{journal}{arXiv preprint arXiv:1504.03594} ,
  \bibinfo{pages}{15}.
\bibitem[{Ghiani et~al.(2004)Ghiani, Laporte and
  Musmanno}]{ghiani2004introduction}
\bibinfo{author}{Ghiani, G.}, \bibinfo{author}{Laporte, G.},
  \bibinfo{author}{Musmanno, R.}, \bibinfo{year}{2004}.
\newblock \bibinfo{title}{Introduction to logistics systems planning and
  control}.
\newblock \bibinfo{publisher}{John Wiley \& Sons}.
\bibitem[{Gilbert(2000)}]{gilbert2000coordination}
\bibinfo{author}{Gilbert, S.M.}, \bibinfo{year}{2000}.
\newblock \bibinfo{title}{Coordination of pricing and multiple-period
  production across multiple constant priced goods}.
\newblock \bibinfo{journal}{Management Science} \bibinfo{volume}{46},
  \bibinfo{pages}{1602--1616}.
\bibitem[{Goyal(1975)}]{goyal1975scheduling}
\bibinfo{author}{Goyal, S.}, \bibinfo{year}{1975}.
\newblock \bibinfo{title}{Scheduling a multi-product single machine system a
  new approach}.
\newblock \bibinfo{journal}{The International Journal of Production Research}
  \bibinfo{volume}{13}, \bibinfo{pages}{487--493}.
\bibitem[{Guastaroba et~al.(2017)Guastaroba, Savelsbergh and
  Speranza}]{guastaroba2017adaptive}
\bibinfo{author}{Guastaroba, G.}, \bibinfo{author}{Savelsbergh, M.},
  \bibinfo{author}{Speranza, M.G.}, \bibinfo{year}{2017}.
\newblock \bibinfo{title}{Adaptive kernel search: A heuristic for solving mixed
  integer linear programs}.
\newblock \bibinfo{journal}{European Journal of Operational Research}
  \bibinfo{volume}{263}, \bibinfo{pages}{789--804}.
\bibitem[{Guastaroba and Speranza(2014)}]{guastaroba2014heuristic}
\bibinfo{author}{Guastaroba, G.}, \bibinfo{author}{Speranza, M.G.},
  \bibinfo{year}{2014}.
\newblock \bibinfo{title}{A heuristic for bilp problems: the single source
  capacitated facility location problem}.
\newblock \bibinfo{journal}{European Journal of Operational Research}
  \bibinfo{volume}{238}, \bibinfo{pages}{438--450}.
\bibitem[{Iyogun(1991)}]{iyogun1991heuristic}
\bibinfo{author}{Iyogun, P.}, \bibinfo{year}{1991}.
\newblock \bibinfo{title}{Heuristic methods for the multi-product dynamic lot
  size problem}.
\newblock \bibinfo{journal}{Journal of the Operational Research Society}
  \bibinfo{volume}{42}, \bibinfo{pages}{889--894}.
\bibitem[{Lee et~al.(2005)Lee, Han and Cho}]{lee2005heuristic}
\bibinfo{author}{Lee, W.S.}, \bibinfo{author}{Han, J.H.}, \bibinfo{author}{Cho,
  S.J.}, \bibinfo{year}{2005}.
\newblock \bibinfo{title}{A heuristic algorithm for a multi-product dynamic
  lot-sizing and shipping problem}.
\newblock \bibinfo{journal}{International Journal of Production Economics}
  \bibinfo{volume}{98}, \bibinfo{pages}{204--214}.
\bibitem[{Linderoth and Savelsbergh(1999)}]{linderoth1999computational}
\bibinfo{author}{Linderoth, J.T.}, \bibinfo{author}{Savelsbergh, M.W.},
  \bibinfo{year}{1999}.
\newblock \bibinfo{title}{A computational study of search strategies for mixed
  integer programming}.
\newblock \bibinfo{journal}{INFORMS Journal on Computing} \bibinfo{volume}{11},
  \bibinfo{pages}{173--187}.
\bibitem[{Lodi and Tramontani(2013)}]{lodi2013performance}
\bibinfo{author}{Lodi, A.}, \bibinfo{author}{Tramontani, A.},
  \bibinfo{year}{2013}.
\newblock \bibinfo{title}{Performance variability in mixed-integer
  programming}, in: \bibinfo{booktitle}{Theory driven by influential
  applications}. \bibinfo{publisher}{INFORMS}, pp. \bibinfo{pages}{1--12}.
\bibitem[{Lu and Qi(2011)}]{lu2011dynamic}
\bibinfo{author}{Lu, L.}, \bibinfo{author}{Qi, X.}, \bibinfo{year}{2011}.
\newblock \bibinfo{title}{Dynamic lot sizing for multiple products with a new
  joint replenishment model}.
\newblock \bibinfo{journal}{European Journal of Operational Research}
  \bibinfo{volume}{212}, \bibinfo{pages}{74--80}.
\bibitem[{Madigan(1968)}]{madigan1968scheduling}
\bibinfo{author}{Madigan, J.}, \bibinfo{year}{1968}.
\newblock \bibinfo{title}{Scheduling a multi-product single machine system for
  an infinite planning period}.
\newblock \bibinfo{journal}{Management Science} \bibinfo{volume}{14},
  \bibinfo{pages}{713--719}.
\bibitem[{McBride and Zufryden(1988)}]{mcbride1988integer}
\bibinfo{author}{McBride, R.D.}, \bibinfo{author}{Zufryden, F.S.},
  \bibinfo{year}{1988}.
\newblock \bibinfo{title}{An integer programming approach to the optimal
  product line selection problem}.
\newblock \bibinfo{journal}{Marketing Science} \bibinfo{volume}{7},
  \bibinfo{pages}{126--140}.
\bibitem[{Meindl and Templ(2012)}]{meindl2012analysis}
\bibinfo{author}{Meindl, B.}, \bibinfo{author}{Templ, M.},
  \bibinfo{year}{2012}.
\newblock \bibinfo{title}{Analysis of commercial and free and open source
  solvers for linear optimization problems}.
\newblock \bibinfo{journal}{Eurostat and Statistics Netherlands within the
  project ESSnet on common tools and harmonised methodology for SDC in the ESS}
  \bibinfo{volume}{20}.
\bibitem[{Minner(2009)}]{minner2009comparison}
\bibinfo{author}{Minner, S.}, \bibinfo{year}{2009}.
\newblock \bibinfo{title}{A comparison of simple heuristics for multi-product
  dynamic demand lot-sizing with limited warehouse capacity}.
\newblock \bibinfo{journal}{International Journal of Production Economics}
  \bibinfo{volume}{118}, \bibinfo{pages}{305--310}.
\bibitem[{Pryor and Chinneck(2011)}]{pryor2011faster}
\bibinfo{author}{Pryor, J.}, \bibinfo{author}{Chinneck, J.W.},
  \bibinfo{year}{2011}.
\newblock \bibinfo{title}{Faster integer-feasibility in mixed-integer linear
  programs by branching to force change}.
\newblock \bibinfo{journal}{Computers \& Operations Research}
  \bibinfo{volume}{38}, \bibinfo{pages}{1143--1152}.
\bibitem[{Salvietti and Smith(2008)}]{salvietti2008profit}
\bibinfo{author}{Salvietti, L.}, \bibinfo{author}{Smith, N.R.},
  \bibinfo{year}{2008}.
\newblock \bibinfo{title}{A profit-maximizing economic lot scheduling problem
  with price optimization}.
\newblock \bibinfo{journal}{European Journal of Operational Research}
  \bibinfo{volume}{184}, \bibinfo{pages}{900--914}.
\bibitem[{Silver(1973)}]{silver1973heuristic}
\bibinfo{author}{Silver, E.A.}, \bibinfo{year}{1973}.
\newblock \bibinfo{title}{A heuristic for selecting lot size quantities for the
  case of a deterministic time-varying demand rate and discrete opportunities
  for replenishment}.
\newblock \bibinfo{journal}{Production Inventory Management}
  \bibinfo{volume}{2}, \bibinfo{pages}{64--74}.
\bibitem[{Silver et~al.(1998)Silver, Pyke, Peterson
  et~al.}]{silver1998inventory}
\bibinfo{author}{Silver, E.A.}, \bibinfo{author}{Pyke, D.F.},
  \bibinfo{author}{Peterson, R.}, et~al., \bibinfo{year}{1998}.
\newblock \bibinfo{title}{Inventory management and production planning and
  scheduling}. volume~\bibinfo{volume}{3}.
\newblock \bibinfo{publisher}{Wiley New York}.
\bibitem[{Sun et~al.(2010)Sun, Huang and Jaruphongsa}]{sun2010economic}
\bibinfo{author}{Sun, H.}, \bibinfo{author}{Huang, H.C.},
  \bibinfo{author}{Jaruphongsa, W.}, \bibinfo{year}{2010}.
\newblock \bibinfo{title}{The economic lot scheduling problem under extended
  basic period and power-of-two policy}.
\newblock \bibinfo{journal}{Optimization Letters} \bibinfo{volume}{4},
  \bibinfo{pages}{157--172}.
\bibitem[{Wagner and Whitin(2004)}]{wagner2004dynamic}
\bibinfo{author}{Wagner, H.M.}, \bibinfo{author}{Whitin, T.M.},
  \bibinfo{year}{2004}.
\newblock \bibinfo{title}{Dynamic version of the economic lot size model}.
\newblock \bibinfo{journal}{Management Science} \bibinfo{volume}{50},
  \bibinfo{pages}{1770--1774}.
\bibitem[{Williams and Redwood(1974)}]{williams1974structured}
\bibinfo{author}{Williams, H.P.}, \bibinfo{author}{Redwood, A.},
  \bibinfo{year}{1974}.
\newblock \bibinfo{title}{A structured linear programming model in the food
  industry}.
\newblock \bibinfo{journal}{Journal of the Operational Research Society}
  \bibinfo{volume}{25}, \bibinfo{pages}{517--527}.
\bibitem[{Witzig and Gleixner(2019)}]{witzig2019conflictdriven}
\bibinfo{author}{Witzig, J.}, \bibinfo{author}{Gleixner, A.},
  \bibinfo{year}{2019}.
\newblock \bibinfo{title}{Conflict-driven heuristics for mixed integer
  programming}.
\newblock \href{http://arxiv.org/abs/1902.02615}{{\tt arXiv:1902.02615}}.

\end{thebibliography}
